\documentclass[a4paper,oneside,10.5pt]{article}%
\usepackage{amsmath}
\usepackage{amsfonts}
\usepackage{amssymb}
\usepackage{graphicx}%
\providecommand{\U}[1]{\protect\rule{.1in}{.1in}}
\newtheorem{theorem}{Theorem}[section]

\newtheorem{definition}[theorem]{Definition}
\newtheorem{assumption}[theorem]{Assumption}

\newtheorem{lemma}[theorem]{Lemma}

\newtheorem{remark}[theorem]{Remark}

\numberwithin{equation}{section}

\begin{document}

\title{Multidimensional viscosity solutions theory of semi-linear partial differential equations}
\author{Shuzhen Yang\thanks{Institute for Financial study, Shandong University, Jinan, Shandong 250100, PR China and Center for Mathematical Economics, Bielefeld University, Germany, (yangsz@sdu.edu.cn). }
\thanks{This work was supported by the Fundamental Research Funds of Shandong University (2015GN015); Supported by the Promotive research fund for excellent young and middle-aged scientists of Shandong Province (BS2015SF004); Supported by the China Postdoctoral Science
Foundation funded project (2015M570584).}}
\date{}
\maketitle

\textbf{Abstract}: In this study, we concern the multidimensional viscosity solutions theory of a kind of  semi-linear partial differential equations (PDEs). A new definition of viscosity solution for this multidimensional semi-linear PDEs which is related to a type of multidimensional backward stochastic differential equations (BSDEs) is given. Further more, we establish the existence and uniqueness results for the viscosity solution of this semi-linear PDEs via the comparison theorem of the related BSDEs and a smooth approximation technique.

\bigskip


\addcontentsline{toc}{section}{\hspace*{1.8em}Abstract}

\section{Introduction}

Since the notion of viscosity solution was invented by Crandall and Lions \cite{CL83} which has become a universal tool to study such a broad fundamental subject. For detailed exposition of such a tool and the related general dynamic programming theory on optimal stochastic control, we refer Crandall, Ishii and Lions \cite{CIL92} for the survey of viscosity solution theory, and the monographs of Fleming and Soner \cite{FS06}; refer \cite{PL79,PL82,PL83a,PL83b} for optimal control of diffusion processes and viscosity solution theory of Hamilton-Jacobi-Bellman equations. For viscosity solution is a big literature, we don't give much more details here.

Ishii and Koike first studied the multidimensional  viscosity solutions of a monotone systems of second-order elliptic partial differential equations (PDEs) in \cite{IK91}. Later Koike showed the uniqueness of viscosity solutions for monotone systems of fully nonlinear PDEs under Dirichlet condition \cite{K94}. In addition, from a stochastic interpretation point of view, Pardoux et al \cite{PPR97}, Buckdahn and Hu \cite{BH10} studied a kind of systems of coupled
Hamilton-Jacobi-Bellman equations and established the related viscosity solution. Also, systems
of coupled Hamilton-Jacobi-Bellman-Isaacs equations and their interpretation through stochastic
differential games, see Qian \cite{Q15}.

In this paper, we aim to study the viscosity solution of the following coupled semi-linear partial differential equation, 
\begin{equation}%
\begin{array}
[c]{l}%
\partial_{t}u^i(t,x)+\mathcal{L}u^i(t,x)
+f(x,u(t,x),\sigma(x)\partial_{x}u^i(t,x))=0;\\
u^i(T,x)=\Phi^{i}(x),\text{ \ \ }x\in\mathbb{R}^d;\text{\ }\\
\end{array}
\label{in-00}%
\end{equation}
where
\[
\mathcal{L}=\frac{1}{2}\sigma\sigma^T(x)\partial_{xx}+b(x)\partial_{x}.
\]

Notice that, the model (\ref{in-00}) is similar with the one in \cite{IK91}, but our assumptions (specially, the conditions (A.1)and (A.2) in \cite{IK91} and Assumption 4.1 in our paper) and method are very different from them. We will deal with our model via a multidimensional BSDEs and a smooth approximation technique.

In order to introduce the multidimensional viscosity  solution theory of semi-linear partial differential equations of our model. Let us quickly scan the theory of backward stochastic differential equations (BSDEs) which be related with semi-linear PDEs, it is well known that the nonlinear BSDEs was first introduced by Pardoux and Peng \cite{PP90}. Independently,
Duffie and Epstein \cite{DE92} presented a stochastic differential recursive
utility which corresponds to the solution of a particular BSDEs. Then the BSDEs
point of view gives a simple formulation of recursive utilities (see
\cite{EPQ97}). Since then, the classical stochastic optimal control problem has been generalized to a so called "stochastic recursive optimal control problem" in which the cost functional is described by the solution of BSDEs. Peng \cite{P92} obtained the Hamilton--Jacobi--Bellman equation for this kind of problem and proved that the value function is its viscosity solution. In \cite{P97}, Peng
generalized his results and originally introduced the notion of stochastic
backward semigroups which allows him to prove the dynamic programming
principle in a very straightforward way. This backward semigroup approach is
proved to be a useful tool for stochastic optimal control problems. Further study see \cite{BL08,BLRT10,WY08}.

Notice that, one dimensional BSDEs is related with the stochastic recursive optimal control problem, and multidimensional BSDEs is related to multidimensional  stochastic recursive problem, i.e., consider the following forward-backward systems (for simplicity, n=2, d=1, more details see Section 2).
\begin{align}
dX_{s}^{t,x}  &  =b(X_{s}^{t,x})ds+\sigma(X_{s}^{t,x})dW_{s}%
,\label{in-1}\\
X_{t}^{t,x}  &  =x,\nonumber
\end{align}
and%
\begin{equation}%
\begin{array}
[c]{ll}%
dY_{s}^{1,t,x}= -f(X_{s}^{t,x},Y_{s}%
^{t,x},Z_{s}^{1,t,x})ds+Z_{s}^{1,t,x}dW_{s};\\
dY_{s}^{2,t,x}= -g(X_{s}^{t,x},Y_{s}%
^{t,x},Z_{s}^{2,t,x})ds+Z_{s}^{2,t,x}dW_{s};\\
Y_{T}^{1,t,x}=\Phi^{1}(X_{T}^{t,x}),\ Y_{T}^{2,t,x}=\Phi^{2}(X_{T}^{t,x}),
\end{array}
\label{in-2}%
\end{equation}
where  $Y_{}%
^{t,x}=(Y_{}%
^{1,t,x},Y_{}%
^{2,t,x})$. Then, we can define the multidimensional utility value functions as,
$$
u(t,x)=Y_{t}^{1,t,x},\ v(t,x)=Y_{t}^{2,t,x}.
$$
If we suppose $(u,v)$ is smooth function, we can relate the value functions $(u,v)$ with
the following partial differential equations,
\begin{equation}%
\begin{array}
[c]{l}%
\partial_{t}u(t,x)+\mathcal{L}u(t,x)
+f(x,u(t,x),v(t,x),\sigma(x)\partial_{x}u(t,x))=0;\\
u(T,x)=\Phi^{1}(x),\text{ \ \ }x\in\mathbb{R};\text{\ }\\
\partial_{t}v(t,x)+\mathcal{L}v(t,x)
+g(x,u(t,x),v(t,x),\sigma(x)\partial_{x}v(t,x))=0;\\
v(T,x)=\Phi^{2}(x),\text{ \ \ }x\in\mathbb{R},
\end{array}
\label{in-3}%
\end{equation}
where
\[
\mathcal{L}=\frac{1}{2}(\sigma(x))^2\partial_{xx}+b(x)\partial_{x}.
\]
Notice that,  $u,v$ are coupled by $f,g$, if we suppose $f,g$ satisfy Lipschatiz conditions (see Section 2, Assumption \ref{ass-2} ), we can't get the relation between the solution of the forward-backward systems (\ref{in-1}), (\ref{in-2}) and the solution of PDEs (\ref{in-3}) by classical solution theory.

In the following, we give a definition of multidimensional viscosity solution for PDEs (\ref{in-3}) which is similar with the classical definition of one dimensional viscosity solution, but the property of which is different from it. Then, we add a  monotonic condition of $f,g$ on $(u,v)$ to guarantee the comparison theorem for BSDEs (\ref{in-2}), by which we prove that the solution of BSDEs (\ref{in-2}) is the viscosity solution of PDEs (\ref{in-3}), this is the existence result. For the uniqueness, we show that the solution of  BSDEs (\ref{in-2}) is the maximum viscosity subsolution and the minimum  viscosity supersolution of PDEs (\ref{in-3}) via another condition on $f,g$ (which is different with the conditions (A1) and (A2) in \cite{IK91}) and a smooth approximation technique.

The paper is organized as follows. In Section 2, we formulate our model and present some fundamental results on BSDEs theory. The existence results for the viscosity solution of PDEs (\ref{in-3}) is established in
Section 3. In Section 4, we prove that the solution of BSDEs (\ref{in-2}) is the unique viscosity solution of PDEs (\ref{in-3}).

\section{Preliminary}
\subsection{Formulation of the problem}

Let $\Omega=C([0,T];\mathbb{R}^{})$ and ${P}$ be the Wiener measure on
$(\Omega,\mathcal{B}(\Omega))$. We denote by $W=(W(t)_{t\in\lbrack0,T]})$ the
canonical Wiener process, with $W(t,\omega)=\omega(t)$, $t\in\lbrack0,T]$,
$\omega\in\Omega$. For any $t\in\lbrack0,T]$, we denote by $\mathcal{F}_{t}$
the ${P}$-completion of $\sigma(W(s),s\in\lbrack0,t])$, and define the following spaces of processes:

$\mathcal{S}^2(0,T;\mathbb{R}^2):=\{\psi \ \text{continuous and progressively measurable}; \ \sup_{0\leq t\leq T}E[\left|\psi(t)\right|^2]< \infty\}$;

$\mathcal{H}^2(0,T;\mathbb{R}):=\{\psi \ \text{ progressively measurable};
\ \ E[\int_0^T\left|\psi(t)\right|^2dt]< \infty\}$;

$L^{2}(\Omega,\mathcal{F}_{T},P):=\{ \xi \ \text{ random variable}; \ \ E[\left|\xi\right|^2]<\infty    \}$.\\
 Let us consider the following forward-backward system. Without loss of generality, we consider the case $n=2$:
\begin{align}
dX_{s}^{t,x}  &  =b(X_{s}^{t,x})ds+\sigma(X_{s}^{t,x})dW_{s}%
;\label{pde-1}\\
X_{t}^{t,x}  &  =x,\nonumber
\end{align}
and%
\begin{equation}%
\begin{array}
[c]{ll}%
dY_{s}^{1,t,x}= -f(X_{s}^{t,x},Y_{s}%
^{t,x},Z_{s}^{1,t,x})ds+Z_{s}^{1,t,x}dW_{s};\\
dY_{s}^{2,t,x}= -g(X_{s}^{t,x},Y_{s}%
^{t,x},Z_{s}^{2,t,x})ds+Z_{s}^{2,t,x}dW_{s};\\
Y_{T}^{1,t,x}=\Phi^{1}(X_{T}^{t,x}),\ Y_{T}^{2,t,x}=\Phi^{2}(X_{T}^{t,x}).
\end{array}
\label{pde-2}%
\end{equation}
where $Y^{t,x}=(Y^{1,t,x},Y^{2,t,x})$ and
\[%
\begin{array}
[c]{l}%
b, \sigma:\mathbb{R}\rightarrow\mathbb{R};\\
f, g:\mathbb{R}\times\mathbb{R}^2\times\mathbb{R}\rightarrow\mathbb{R}^{};\\
\Phi^{1},\Phi^{2}:\mathbb{R}\rightarrow\mathbb{R}.
\end{array}
\]

 In order to obtain the well-posedness of Forward-Backward SDEs (\ref{pde-1}) and (\ref{pde-2}), we denote by $h:=f,g;\ \Phi:=\Phi^1,\Phi^2$ and assume that $b,\sigma,h,\Phi$ are deterministic functions and satisfy the
following conditions:
\begin{assumption}
\label{ass-1}There exists a constant $C>0$ such that%
\[%
\begin{array}
[c]{l}%
\mid b(x^{1})-b(x^{2})\mid+\mid\sigma(x^{1})-\sigma(x^{2})\mid\leq
C\mid x^{1}-x^{2}\mid;\\
\mid\Phi(x^{1})-\Phi(x^{2})\mid\leq C\mid x^{1}-x^{2}\mid,
\end{array}
\]
$\forall x^{1},x^{2}\in\mathbb{R}$.
\end{assumption}

\begin{remark}
Suppose $b,\sigma$ satisfy Assumption \ref{ass-1}. Then there exists a constant $C>0$
such that%
\[%
\begin{array}
[c]{l}%
\mid b(x)\mid+\mid\sigma(x)\mid\leq C(1+\mid x\mid);\\
\mid\Phi(x)\mid\leq C(1+\mid x\mid),\text{ \ \ }\forall x\in\mathbb{R}.
\end{array}
\]

\end{remark}
\begin{assumption}
\label{ass-2}Let $h$ be Lipschitz continuous in $x,y,z$, i.e.%
\[
\left\vert h(x_1,y_1,z_{1})-h(x_2,y_2,z_{2})\right\vert \leq C\big{(}\left|x_{1}-x_{2}\right|+\left|y_{1}-y_{2}\right|+\left|z_{1}-z_{2}\right|\big{)},
\]
and linear growth, i.e. there exists a positive constant $C_{}$ such that
\[
\left\vert h(x_1,y_1,z_{1})\right\vert
\leq C_{}(\left\vert x_1\right\vert+\left\vert y_1\right\vert+\left\vert z_{1}\right\vert +1),
\]
$\forall \ x=(x_1,x_2), y_1,y_2,z=(z_1,z_2)\in\mathbb{R}^{2}.$
\end{assumption}

The following assumption guarantees the comparison theorem of BSDEs (\ref{pde-2}).
\begin{assumption}
\label{ass-3}Let $f$ and $g$ satisfy the following monotonic conditions, i.e.,
\[
\begin{array}
[c]{cl}
f(x,y,y_{1},z)\geq f(x,y,y_{2},z);\\
g(x,y_{1},y,z)\geq g(x,y_{2},y,z),\\
\end{array}
\]
$\forall  \ x, y_1,y_2,y,z\in\mathbb{R}$ with $ y_1\geq y_2$.
\end{assumption}

Under the above assumptions, we have the following theorems.
\begin{theorem}
Let $b,\sigma$ satisfy Assumption \ref{ass-1}, then there exists a unique adapted solution $X$ for equation (\ref{pde-1}).
\end{theorem}
Based on the theory of BSDEs, we have the next existence and uniqueness results. We refer \cite{PP90} and \cite{EPQ97} for the theory of BSDEs.
\begin{theorem}
\label{th-cu} Let $b,\sigma$ satisfy Assumption \ref{ass-1} and $f,g$\ satisfy Assumption \ref{ass-2}, then there exists a unique adapted solutions
$(Y,Z)$ in $\mathcal{S}^2\times\mathcal{H}^2$ solving equation (\ref{pde-2}).
\end{theorem}

The following comparison theorem for equations (\ref{cp-11}) and (\ref{cp-12}) is the key step when proving the existence of the viscosity solution of semi-linear PDEs. For reader convenience, we show the proof in the Appendix. %
\begin{equation}%
Y_{1,i}(t)=  \xi_{1,i}+\int_{t}^{T}f_{i}(Y_{i}(s),Z_{1,i}(s))ds-\int%
_{t}^{T}Z_{1,i}(s)dW(s);
\label{cp-11}%
\end{equation}
\begin{equation}
Y_{2,i}(t)=  \xi_{2,i}+\int_{t}^{T}g_{i}(Y_{i}(s),Z_{2,i}(s))ds-\int%
_{t}^{T}Z_{2,i}(s)dW(s),
\label{cp-12}%
\end{equation}
where $Y_i=(Y_{1,i},Y_{2,i})$, $i=1,2$.

\begin{theorem}
\label{theorem-cp}If $f_i,g_i$ satisfy Assumptions \ref{ass-2} and \ref{ass-3}; $\xi_{1,i}%
,\xi_{2,i}\in L^{2}(\Omega,\mathcal{F}_{T},P)$ and
\begin{equation}%
\begin{array}
[c]{rl}%
f_1\leq f_2,g_1 \leq g_2;\quad
\xi_{1,1}\leq \xi_{1,2},\xi_{2,1}\leq \xi_{2,2}, \quad \mathbb{P}-a.s.
\end{array}
\label{cp-3}%
\end{equation}
then we have comparison results
\[
{Y}_{1,1}\leq {Y}_{1,2},\quad {Y}_{2,1}\leq {Y}_{2,2},\quad a.s.,a.e.
\]
\end{theorem}

\subsection{Classical solution of Forward Backward SDEs}

Let $t\in\lbrack0,T]$ and $x\in\mathbb{R}$. Given the following forward-backward SDEs:
\begin{equation}
\left \{
\begin{array}
[c]{ll}
X^{t,x}_s=x+\int_{t}^{s}b(X^{{t}%
,x}_r)dr+\int_{t}^{s}\sigma
(X^{{t},x}_r)dW_r;\\
Y^{{t},x}_s=g(X^{{t},x}_T)-\int_{s}%
^{T}h(X^{{t},x}_r,Y^{{t},x}_r,Z^{{t},x}_r)dr-\int_{s}^{T}Z^{{t},x}_rdW_r,
\label{2.1}%
\end{array}
\right.
\end{equation}%
for every $s\in\lbrack t,T]$.

We  recall some notions and results in Pardoux and Peng \cite{PP92}.
$\mathcal{C}^{n}(\mathbb{R};\mathbb{R})$, $\mathcal{C}_{b}^{n}(\mathbb{R}; \mathbb{R})$, $\mathcal{C}_{p}^{n}(\mathbb{R};\mathbb{R})$ will denote
respectively the set of functions of class $\mathcal{C}^{n}$ from $\mathbb{R}$ into
$\mathbb{R}$, the set of those functions of class $\mathcal{C}_{b}^{n}$ whose
partial derivatives of order less than or equal to $n$ are bounded, and the
set of those functions of class $\mathcal{C}_{p}^{n}$ which, together with all their
partial derivatives of order less than or equal to $n$, grow at most like a
polynomial function of the variable $x$ at infinity.

\begin{assumption}\label{ass1}
$b,\sigma \in \mathcal{C}_{l,b}^{3}(\mathbb{R};\mathbb{R})$ with the first order partial derivatives in x are bounded, as well
as their derivatives of order one and two with respect to x.
\end{assumption}

\begin{assumption}\label{ass2}
 $h,g$ are of class $\mathcal{C}_{p}^{3}$ and the first order partial
derivatives in $x,y$ and $z$ are bounded, as well as their derivatives of up to
order two with respect to $x,y,z$.
\end{assumption}

In Pardoux and Peng \cite{PP92}, under Assumptions \ref{ass1} and \ref{ass2},
$Y$ and $Z$ are related in the following sense:
\begin{equation*}
\left\{
\begin{array}
[c]{ll}
Y^{t,x}_s=u(s,X_s^{t,x});\\
Z^{{t},x}_s=\partial_{x}u(s,X^{{t},x}_s)\sigma(X^{{t},x}_s),\quad P-a.s.
\end{array}
\right.
\end{equation*}
Equation (\ref{2.1}) be related to the following partial differential equation:
\begin{equation}%
\left\{
\begin{array}
[c]{ll}%
\partial_{t}u({t},x)+\mathcal{L}u({t},x)=h(x,u({t},x), \partial_{x}u({t},x)\sigma(x));\\
u({T},x)=g(x),
\end{array}
\right.
\label{2.3}%
\end{equation}
where $\partial$ is the gradient operator and
\[
\mathcal{L}=\frac{1}{2}(\sigma(x))^2\partial_{xx}^2+b(x)\partial_{x}.
\]
\begin{theorem}
Suppose Assumptions \ref{ass1} and \ref{ass2} hold. If $u$ belongs to
$\mathcal{C}^{1,2}$ and $(u,v)$ is the solution of
equation (\ref{2.3}) such that $(u,v)$ is uniformly Lipschitz continuous and
bounded by $C(1+|x|)$, then we have $u({t},x)=Y^{{t},x}(t)$, for each $({t},x)\in[0,T]
\times\mathbb{R}$, where $(X^{{t},x}(s),Y^{{t},x}(s),Z^{{t},x}(s))_{t\leq s\leq T}$ is the unique solution of forward Backward SDE (\ref{2.1}).
\end{theorem}

\begin{theorem}
\label{theorem-cs}
Under Assumptions \ref{ass1} and \ref{ass2}, the function $u({t},x)=Y^{{t},x}(t)$ is the unique ${\mathcal{C}}^{1,2}$-solution of PDE
(\ref{2.3}).
\end{theorem}

\section{The existence of viscosity solution}
In this section, we will prove that the solution of BSDEs solves the multidimensional PDEs in the mean of viscosity solution.
Consider the following PDEs:%
\begin{equation}%
\begin{array}
[c]{l}%
\partial_{t}u(t,x)+\mathcal{L}u(t,x)+f(x,u,v,\sigma(x)\partial_{x}u(t,x))=0;\\
u(T,x)=\Phi^{1}(x),\text{ \ \ }x\in\mathbb{R};\text{\ }\\
\partial_{t}v(t,x)+\mathcal{L}v(t,x)+g(x,u,v,\sigma(x)\partial_{x}v(t,x))=0;\\
v(T,x)=\Phi^{2}(x),\text{ \ \ }x\in\mathbb{R},
\end{array}
\label{vis-1}%
\end{equation}
where
\[
\mathcal{L}=\frac{1}{2}(\sigma(x))^2\partial_{xx}+b(x)\partial_{x}.
\]

For given initial data $(t,x)$, we define
\[
u(t,x):=Y_{t}^{1,t,x},v(t,x):=Y_{t}^{2,t,x},
\]
 where $(Y^{1,t,x},Y^{2,t,x})$ is the solution of BSDEs (\ref{pde-2}). By the markov property of BSDEs (\ref{pde-1}) and (\ref{pde-2}), we have the following results:
\begin{lemma}
\label{dpp-1}
For $t\leq s$, we have
\[
u(s,X_s^{t,x})=Y_{s}^{1,t,x},\ v(s,X_s^{t,x})=Y_{s}^{2,t,x}.
\]
\end{lemma}
\textbf{Proof:} This lemma really come from the paper \cite{P97}, so we omit it.

\bigskip

Now, we first give the definition of the viscosity solution of PDEs
(\ref{vis-1}) as follows.
\begin{definition}
\label{de-1}
Let $w=(w_1,w_2)\in\mathcal{C}([0,T]\times \mathbb{R})\times \mathcal{C}([0,T]\times \mathbb{R})$, we say that $w$ is a viscosity subsolution of
(\ref{vis-1}), if $\forall \ \Gamma=(\Gamma_1,\Gamma_2)\in\mathcal{C}^{1,2}([0,T]\times \mathbb{R})\times \mathcal{C}^{1,2}([0,T]\times \mathbb{R})$ at $(t_1,x_1),(t_2,x_2)$ satisfying $\Gamma_1\geq w_1$ and $\Gamma_2\geq w_2$ on $[0,T]\times \mathbb{R}$ and $\Gamma_1(t_1,x_1)=w_1(t_1,x_1),\  \Gamma_2(t_2,x_2)=w_2(t_2,x_2)$, we have
\begin{equation*}
\begin{array}
[c]{ll}
\partial_{t}\Gamma_1(t_1,x_1)+\mathcal{L}\Gamma_1(t_1,x_1)+f(x_1,\Gamma_1,\Gamma_2,\sigma(x_1)\partial_{x}\Gamma_1(t_1,x_1))\geq0;\\
\partial_{t}\Gamma_2(t_2,x_2)+\mathcal{L}\Gamma_2(t_2,x_2)+g(x_2,\Gamma_1,\Gamma_2,\sigma(x_2)\partial_{x}\Gamma_2(t_2,x_2))\geq0.\\
\end{array}
\end{equation*}

We say that $w$ is a viscosity supersolution of (\ref{vis-1}), if $\forall \ \Gamma=(\Gamma_1,\Gamma_2)\in\mathcal{C}^{1,2}([0,T]\times \mathbb{R})\times \mathcal{C}^{1,2}([0,T]\times \mathbb{R})$ at $(t_1,x_1),(t_2,x_2)$  satisfying $\Gamma_1\leq w_1$ and $\Gamma_2\leq w_2$ on $[0,T]\times \mathbb{R}$ and $\Gamma_1(t_1,x_1)=w_1(t_1,x_1),\  \Gamma_1(t_2,x_2)=w_2(t_2,x_2)$, we have
\begin{equation*}
\begin{array}
[c]{ll}
\partial_{t}\Gamma_1(t_1,x_1)+\mathcal{L}\Gamma_1(t_1,x_1)+f(x_1,\Gamma_1,\Gamma_2,\sigma(x_1)\partial_{x}\Gamma_1(t_1,x_1))\leq0;\\
\partial_{t}\Gamma_2(t_2,x_2)+\mathcal{L}\Gamma_2(t_2,x_2)+g(x_2,\Gamma_1,\Gamma_2,\sigma(x_2)\partial_{x}\Gamma_2(t_2,x_2))\leq0.\\
\end{array}
\end{equation*}

We say that $w$ is a viscosity solution of (\ref{vis-1}) if it is both a
viscosity supersolution and a viscosity subsolution of (\ref{vis-1}).
\end{definition}

\begin{remark}
Notice that, the Definition \ref{de-1} for multidimensional viscosity solution  is different from one dimensional case for $\Gamma_1-w_1$ and $\Gamma_1-w_2$ may take minimum (maximum) values at different points, i.e., $(t_1,x_1)\neq (t_2,x_2)$. Thus, we may consider the following definition for viscosity solution:

Let $w=(w_1,w_2)\in\mathcal{C}([0,T]\times \mathbb{R})\times \mathcal{C}([0,T]\times \mathbb{R})$, we say that $w$ is a viscosity subsolution of
(\ref{vis-1}), if $\forall \ \Gamma=(\Gamma_1,\Gamma_2)\in\mathcal{C}^{1,2}([0,T]\times \mathbb{R})\times \mathcal{C}^{1,2}([0,T]\times \mathbb{R})$ at $(t_1,x_1),(t_2,x_2)$ satisfying $\Gamma_1\geq w_1$ and $\Gamma_2\geq w_2$ on $[0,T]\times \mathbb{R}$ and $\Gamma_1(t_1,x_1)=w_1(t_1,x_1),\  \Gamma_1(t_2,x_2)=w_2(t_2,x_2)$, we have
\begin{equation*}
\begin{array}
[c]{ll}
\partial_{t}\Gamma_1(t_1,x_1)+\mathcal{L}\Gamma_1(t_1,x_1)+f(x_1,w_1,w_2,\sigma(t_1,x_1)\partial_{x}\Gamma_1(t_1,x_1))\geq0;\\
\partial_{t}\Gamma_2(t_2,x_2)+\mathcal{L}\Gamma_2(t_2,x_2)+g(x_2,w_1,w_2,\sigma(t_2,x_2)\partial_{x}\Gamma_2(t_2,x_2))\geq0.\\
\end{array}
\end{equation*}
Notice that, $\Gamma_1-w_1$ and $\Gamma_1-w_2$ may take minimum values at different points, thus the above definition of multidimensional viscosity solution is different from the one in Definition \ref{de-1}. In this study, we does not pay attention to the property of definition of viscosity solution and we will choose the Definition \ref{de-1} to study the viscosity solution of PDEs (\ref{vis-1}).
\end{remark}

Before to show the solution of BSDEs (\ref{pde-2}) is the multidimensional viscosity solution of the semi-nonlinear PDEs (\ref{vis-1}), let us give some lemmas which are used later.

Given  $\Gamma_1,\Gamma_2\in \mathcal{C}^{1,2}([t,T]\times \mathbb{R})$ and assume $f,g$ satisfy the Lipschitz and linear growth conditions.
\begin{equation*}
\begin{array}
[c]{ll}
F_{11}(s,x,y_1,y_2,z_1)=&\partial_s\Gamma_1(s,x)+\mathcal{L}\Gamma_1(s,x)\\
&+f(x,y_1+\Gamma_1(s,x),y_2+\Gamma_2(s,x),z_1+\partial_x\Gamma_1(s,x)\sigma(x));\\
F_{12}(s,x,y_1,y_2,z_2)=&\partial_s\Gamma_2(s,x)+\mathcal{L}\Gamma_2(s,x)\\
&+g(x,y_1+\Gamma_1(s,x),y_2+\Gamma_2(s,x),z_2+\partial_x\Gamma_2(s,x)\sigma(x)),\\
\end{array}
\end{equation*}
$\forall \ (s,x,y_1,y_2,z_1,z_2)\in[t,T]\times\mathbb{R}\times\mathbb{R}
\times\mathbb{R}\times\mathbb{R}\times\mathbb{R}$.

Consider the following BSDEs: $s\in [t,t+\delta]$,
\begin{equation}
\begin{array}
[c]{ll}
Y^1_s=\int_s^{t+\delta}F_{11}(r,X_r^{t,x},Y^1_r,Y^2_r,Z^1_r)dr-\int_s^{t+\delta}Z^1_rdW_r;\\
Y^2_s=\int_s^{t+\delta}F_{12}(r,X_r^{t,x},Y^1_r,Y^2_r,Z^2_r)dr-\int_s^{t+\delta}Z^2_rdW_r,\\
\end{array}
\end{equation}
and
\begin{equation}
\begin{array}
[c]{ll}
Y_s^{1,0}=\Gamma_1(t+\delta,X^{t,x}_{t+\delta})
+\int_s^{t+\delta}f(X_r^{t,x},Y^{1,0}_r,Y^{2,0}_r,Z^{1,0}_r)dr-\int_s^{t+\delta}Z^{1,0}_rdW_r;\\
Y_s^{2,0}=\Gamma_2(t+\delta,X^{t,x}_{t+\delta})
+\int_s^{t+\delta}g(X_r^{t,x},Y^{1,0}_r,Y^{2,0}_r,Z^{2,0}_r)dr-\int_s^{t+\delta}Z^{2,0}_rdW_r.\\
\end{array}
\label{eq-0}
\end{equation}
Applying It\^{o} formula, we obtain the following lemma:
\begin{lemma}
\label{le-1}
For any $s\in[t,t+\delta]$, we have
\begin{equation}
\begin{array}
[c]{ll}
Y^1_s=Y_s^{1,0}-\Gamma_1(s,X_s^{t,x});\\
Y^2_s=Y_s^{2,0}-\Gamma_2(s,X_s^{t,x}).\\
\end{array}
\end{equation}
\end{lemma}
Then, consider the BSDEs: $s\in [t,t+\delta]$
 \begin{equation}
\begin{array}
[c]{ll}
Y^{1,1}_s=\int_s^{t+\delta}F_{11}(r,x,Y^{1,1}_r,Y^{2,1}_r,Z^{1,1}_r)dr-\int_s^{t+\delta}Z^{1,1}_rdW_r;\\
Y^{2,1}_s=\int_s^{t+\delta}F_{12}(r,x,Y^{1,1}_r,Y^{2,1}_r,Z^{2,1}_r)dr-\int_s^{t+\delta}Z^{2,1}_rdW_r.\\
\end{array}
\end{equation}
By basic theory of SDE and BSDEs, we can obtain the next results:
\begin{lemma}
\label{vl-1}
We have
$$
\mid Y_t^{1}-Y^{1,1}_t \mid \leq C\delta^{\frac{3}{2}},\quad \mid Y_t^{2}-Y^{2,1}_t \mid \leq C\delta^{\frac{3}{2}}.
$$
\end{lemma}
Also, we have the representation of $Y^{1,1}$ and $Y^{2,1}$.
\begin{lemma}
\label{vl-2}
We have
$$
Y_t^{1,1}=Y_t^{0,1},\quad  Y_t^{2,1}=Y_t^{0,2},
$$
where $Y^{0,1}$ and $Y^{0,2}$ are the solution of the following differential equations
 \begin{equation}
\begin{array}
[c]{ll}
-dY^{0,1}_s=\int_s^{t+\delta}F_{11}(s,x,Y^{0,1},Y^{0,2},0,0)dr,\quad Y^{0,1}_{t+\delta}=0,\quad s\in[t,t+\delta];\\
-dY^{0,2}_s=\int_s^{t+\delta}F_{12}(s,x,Y^{0,1},Y^{0,2},0,0)dr,\quad Y^{0,2}_{t+\delta}=0,\quad s\in[t,t+\delta],\\
\end{array}
\end{equation}
\end{lemma}
\textbf{Proof:} Notice that, $F_{11}$ and $F_{12}$ are deterministic functions, we obtain that $Z_s^{1,1}=0$, $Z_s^{2,1}=0$ and $Y_s^{1,1}=Y_s^{0,1}$, $Y_s^{2,1}=Y_s^{0,2}$, $s\in[t,t+\delta]$.

By the basic theory of BSDEs, we have the following estimates of the continuity of  $(u(t,x),v(t,x))$ with respect to $t,x$, see Peng {\cite{P92,P97}}.
\begin{lemma}
\label{Lem1} There exists a constant $C_{}>0$ such that, $\forall
t\in\lbrack0,T]$ and $x,x^{\prime}\in\mathbb{R}^{}{,}$ we have

(i) $\mid u(t,x)-u(t,x^{\prime})\mid+\mid v(t,x)-v(t,x^{\prime})\mid\leq C_{}\mid x-x^{\prime}\mid;$

(ii) $\mid u(t,x)\mid+\mid v(t,x)\mid\leq C_{}(1+\mid x\mid).$
\end{lemma}
\begin{lemma}
\label{Lem2}The function $(u,v)$ is $\frac{1}{2}$ H\"{o}lder
continuous in $t$.
\end{lemma}
\noindent\textbf{Proof. }Set $(t,x)\in\mathbb{R}^{}\times\lbrack0,T]$ and $\delta>0.$
By Lemma \ref{dpp-1}, $u(s,X_s^{t,x})=Y_s^{1,t,x}$ and $ v(s,X_s^{t,x})=Y_t^{2,t,x}$, $\forall$ $\varepsilon>0,$ we have
\begin{equation}%
\begin{array}
[c]{ll}%
u(t,x)=u(t+\delta,X_{t+\delta}^{t,x}) +\int_t^{t+\delta}f(X_{s}^{t,x},Y_{s}%
^{t,x},Z_{s}^{1,t,x})ds-\int_t^{t+\delta}Z_{s}^{1,t,x}dW_{s};\\
v(t,x)=v(t+\delta,X_{t+\delta}^{t,x}) +\int_t^{t+\delta}g(X_{s}^{t,x},Y_{s}%
^{t,x},Z_{s}^{2,t,x})ds-\int_t^{t+\delta}Z_{s}^{2,t,x}dW_{s}.\\
\end{array}
\label{ct-1}
\end{equation}
We first show that there exists $C>0$ such that $\mid u(t+\delta,x)-u(t,x)\mid\leq
C\delta^{\frac{1}{2}}.$
By the first equality of (\ref{ct-1}), we have%
\begin{equation}%
\begin{array}
[c]{c}%
u(t+\delta,x)-u(t,x)= I_{\delta}^{1}+I_{\delta}^{2},
\end{array}
\label{ct-2}%
\end{equation}
where
\[%
\begin{array}
[c]{cl}%
I_{\delta}^{1}= & u(t+\delta,x)-u(t+\delta,X_{t+\delta}^{t,x}),\\
I_{\delta}^{2}= &-\int_t^{t+\delta}f(X_{s}^{t,x},Y_{s}%
^{t,x},Z_{s}^{1,t,x})ds+\int_t^{t+\delta}Z_{s}^{1,t,x}dW_{s}.
\end{array}
\]
By Lemma \ref{Lem1}, note that $u$ is 1-H\"{o}lder continuous in
$x$. We have%
\[
{{E}}\big{[}\left\vert I_{\delta}^{1}\right\vert  \big{]}
\leq{{E}}\big{[}\left\vert
u(t+\delta,x)-u(t+\delta,X_{t+\delta}^{t,x,u})\right\vert
\big{]}\leq C{{E}}\big{[}\left\vert X_{t+\delta}^{t,x,u}-x\right\vert\big{]}.
\]
Then by ${{E}}\big{[}\left\vert X_{t+\delta}^{t,x,u}-x\right\vert ^{2}\big{]}\leq
C\delta$ ($C$ will change line by line),%
\[
\left\vert I_{\delta}^{1}\right\vert \leq C\delta^{\frac{1}{2}}.
\]
According to BSDEs (\ref{ct-1}), $I_{\delta
}^{2}$ can be rewritten as
\[%
\begin{array}
[c]{cl}%
{{E}}\big{[} \left\vert I_{\delta}^{2} \right\vert\big{]}
= {{E}}\big{[}\left\vert\int_t^{t+\delta}f(X_{s}^{t,x},Y_{s}%
^{t,x},Z_{s}^{1,t,x})ds-\int_t^{t+\delta}Z_{s}^{1,t,x}dW_{s}\right\vert\big{]}.
\end{array}
\]
It yields that%
\[%
\begin{array}
[c]{rl}%
\left\vert I_{\delta}^{2}\right\vert \leq & \delta^{\frac{1}{2}}%
\bigg{\{}\big{[}{{E}}\int_{t}^{t+\delta}\left\vert f(X_{s}^{t,x,u}%
,Y_{s}^{t,x,u},Z_{s}^{1,t,x})\right\vert ^{2}ds\big{]}^{\frac{1}{2}}+\big{[}{{E}}\int_{t}^{t+\delta}\left\vert Z_{s}^{1,t,x}\right\vert ^{2}ds\big{]}^{\frac{1}{2}}\bigg{\}}
\leq  C\delta^{\frac{1}{2}}.
\end{array}
\]
Thus, we have
\[
\left\vert u(t+\delta,x)-u(t,x) \right\vert \leq C\delta^{\frac{1}{2}}.
\]
Similarly, we can prove that
\[
\left\vert v(t+\delta,x)-v(t,x) \right\vert \leq C\delta^{\frac{1}{2}}.
\]

This completes the proof. $\ \ \ \Box$

\begin{theorem}
\label{viscosity} Let $b,\sigma$ satisfy  Assumptions \ref{ass-1}; $f,g$ satisfy  Assumptions \ref{ass-2} and \ref{ass-3}. Then $(u,v)$ is viscosity
solution of the semi-linear partial differential equation (\ref{vis-1}).
\end{theorem}
\textbf{Proof:} By Lemma \ref{Lem2}, we have that $(u,v)$ is a continuous function of $(t,x)\in [0,T]\times \mathbb{R}$. Next, we prove that $(u,v)$ is viscosity solution of the PDEs (\ref{vis-1}).

\textbf{Step 1}: For any $\Gamma=(\Gamma_1,\Gamma_2)\in\mathcal{C}^{1,2}%
([0,T]\times \mathbb{R})\times \mathcal{C}^{1,2}%
([0,T]\times \mathbb{R})$, let $\Gamma_1\geq u$ and $\Gamma_2\geq v$ on $[0,T]\times \mathbb{R}$ and $(\Gamma_1(t,x),\Gamma_2(t_0,x_0))=(u(t,x),v(t_0,x_0))$ with $t\leq t_0$.

By equation (\ref{vis-1}) and Lemma \ref{dpp-1}, we have
\begin{equation}
\begin{array}
[c]{ll}
u(t,x)= E[\int_{t}^{t+\delta}f(X_{r}^{t,x},Y_r^{1,t,x},Y_r^{2,t,x},Z_r^{1,t,x})dr+u(t+\delta,X_{t+\delta}%
^{t,x})];\\
v(t,x)= E[\int_{t}^{t+\delta}g(X_{r}^{t,x},Y_r^{1,t,x},Y_r^{2,t,x},Z_r^{2,t,x})dr+v(t+\delta,X_{t+\delta}%
^{t,x})].\\
\end{array}
\label{Ito formula1}%
\end{equation}

Notice that $\Gamma_1\geq u,\Gamma_2\geq v$  and $(\Gamma_1(t,x),\Gamma_2(t_0,x_0))=(u(t,x),v(t_0,x_0))$, by comparison Theorem \ref{theorem-cp}, we can compare equations (\ref{eq-0}) and (\ref{Ito formula1}), thus
$$
Y_t^{1,0}\geq u(t,x)=\Gamma_1(t,x).
$$
Then, by Lemma \ref{le-1}, we have
$$
Y_t^{1}=Y_t^{1,0}-\Gamma_1(t,x)\geq 0.
$$
By Lemma \ref{vl-1} and Lemma \ref{vl-2}, we get
$$
Y_t^{1,1}\geq -C\delta^{\frac{3}{2}}.
$$
and
$$
Y_t^{0,1}\geq -C\delta^{\frac{3}{2}}.
$$
Thus
\begin{equation}
\begin{array}
[c]{ll}
-C\delta^{\frac{1}{2}}\leq \delta^{-1}Y^{0,1}_t=\delta^{-1}\int_t^{t+\delta}F_{11}(r,x,Y^{0,1}_r,Y^{0,2}_r,0)dr;\\
\end{array}
\end{equation}
Letting $\delta \to 0$, we get
\begin{equation}
\begin{array}
[c]{ll}
\partial_{t}\Gamma_1(t,x)+\mathcal{L}\Gamma_1(t,x)+f(x,\Gamma_1(t,x),\Gamma_2(t,x),\sigma(x)\partial_x\Gamma_1(t,x))
\geq0.\\
\end{array}
\end{equation}

Similarly, using the same method as above, we consider the following BSDEs:
\begin{equation}
\begin{array}
[c]{ll}
u(t_0,x_0)= E[\int_{t_0}^{t_0+\delta}f(X_{r}^{t_0,x_0},Y_r^{1,t_0,x_0},
Y_r^{2,t_0,x_0},Z_r^{1,t_0,x_0})dr+u(t_0+\delta,X_{t_0+\delta}%
^{t_0,x_0})];\\
v(t_0,x_0)= E[\int_{t_0}^{t+\delta}g(X_{r}^{t_0,x_0},Y_r^{1,t_0,x_0},
Y_r^{2,t_0,x_0},Z_r^{2,t_0,x_0})dr+v(t_0+\delta,X_{t_0+\delta}%
^{t_0,x_0})],\\
\end{array}
\end{equation}
 and which deduce that
\begin{equation}
\begin{array}
[c]{ll}
-C\delta^{\frac{1}{2}}\leq \delta^{-1}Y^{0,2}_{t_0}=\delta^{-1}\int_{t_0}^{t_0+
\delta}F_{12}(r,x_0,Y^{0,1}_r,Y^{0,2}_r,0)dr;\\
\end{array}
\end{equation}
Letting $\delta \to 0$, we have
\begin{equation}
\begin{array}
[c]{ll}
&\partial_{t}\Gamma_2(t_0,x_0)+\mathcal{L}\Gamma_2(t_0,x_0)
+g(x_0,\Gamma_1(t_0,x_0),\Gamma_2(t_0,x_0),\sigma(x_0)\partial_x\Gamma_2(t_0,x_0))
\leq0.\\
\end{array}
\end{equation}
Thus, we complete the proof for this assertion.

\textbf{Step 2}: Next, using the same method as above, we can prove $(u,v)$ is the viscosity supersolution of equation (\ref{vis-1}). Thus, $(u,v)$ is viscosity solution of equation (\ref{vis-1}).

This completes the proof. $\ \ \ \ \ \ \ \ \  \Box$

\section{Uniqueness of viscosity solution}
Based on the results in Section 3, we are now at the stage to prove that the solutions of BSDEs (\ref{pde-2}) is the unique viscosity solution of  the multidimensional PDEs (\ref{vis-1}). Following the Definition \ref{de-1}, once we prove that the solutions of BSDEs (\ref{pde-2}) is the maximum visocisty subsolution and minimum viscosity supersolution of PDEs (\ref{vis-1}), then, any viscosity solution of PDEs (\ref{vis-1}) must be the solution of BSDEs (\ref{pde-2}). In addition, we need the following assumption for proving the uniqueness results.
\begin{assumption}
\label{ass-4}Let $f$ and $g$ satisfy the following monotonic conditions, and there exist constants $C_2>C_1>0$ such that
\[
\begin{array}
[c]{cl}
f(x,y_{11},y_{12},z)-f(x,y_{21},y_{22},z)\leq C_2(y_{21}-y_{11})+C_1(y_{12}-y_{22});\\
g(x,y_{11},y_{12},z)-g(x,y_{21},y_{22},z)\leq C_2(y_{22}-y_{12})+C_1(y_{11}-y_{21}),\\
\end{array}
\]
with $ y_1=(y_{11},y_{12}),y_2=(y_{21},y_{22}),\forall (x,y_1,y_2,z)\in\mathbb{R}\times\mathbb{R}^2\times\mathbb{R}^2\times\mathbb{R}$.
\end{assumption}

Before proving the main results of this section, let us show some preliminary results. Firstly, we will construct some smooth functions which are used to approximate $b,\sigma,f,g,\Phi$. We denote
$$
\psi_{\varepsilon}(y)=\frac{1}{\sqrt{2\pi \varepsilon^2}}e^{-\frac{y^2}{2\varepsilon^2}},\ y\in \mathbb{R},
$$
and the convolution of $b,\sigma,f,g,\Phi^1,\Phi^2$ by
\begin{equation}
\label{unae-1}
\begin{array}
[c]{rl}
{b}_{\varepsilon}(x)=&\int_{\mathbb{R}}b(x^0)\psi_{\varepsilon}(x^0-x)dx^0;\\
{\sigma}_{\varepsilon}(x)=&\int_{\mathbb{R}}\sigma(x^0)\psi_{\varepsilon}(x^0-x)dx^0;\\
{\Phi}^1_{\varepsilon}(x)=&\int_{\mathbb{R}}{\Phi}^1(x^0)\psi_{\varepsilon}(x^0-x)dx^0;\\
{\Phi}^2_{\varepsilon}(x)=&\int_{\mathbb{R}}{\Phi}^2(x^0)\psi_{\varepsilon}(x^0-x)dx^0;\\
{f}_{\varepsilon}(x,y_1,y_2,z)=&\int_{\mathbb{R}^4}\big{[}f(x,y_1,y_2,z)
\psi_{\varepsilon}(x^0-x)\\
&\times\psi_{\varepsilon}(y_1^0-y_1)
\psi_{\varepsilon}(y_2^0-y_2)\psi_{\varepsilon}(z^0-z)\big{]}dx^0dy_1^0dy_2^0dz^0;\\
{g}_{\varepsilon}(x,y_1,y_2,z)=&\int_{\mathbb{R}^4}\big{[}g(x,y_1,y_2,z)
\psi_{\varepsilon}(x^0-x)\\
&\times\psi_{\varepsilon}(y_1^0-y_1)
\psi_{\varepsilon}(y_2^0-y_2)\psi_{\varepsilon}(z^0-z)\big{]}dx^0dy_1^0dy_2^0dz^0,\\
\end{array}
\end{equation}
with $(x,y_1,y_2,z)\in \mathbb{R}\times\mathbb{R}\times\mathbb{R}\times\mathbb{R}$. Thus, we have the following lemma.
\begin{lemma}
\label{visle-1}
There is constant $C>0$ such that
\begin{equation*}
\begin{array}
[c]{ll}
|{b}_{\varepsilon}(x)-b(x)|
+\left|{\sigma}_{\varepsilon}(x)-\sigma(x)\right|\leq C\varepsilon;\\
|{\Phi}^1_{\varepsilon}(x)-\Phi^1(x)|+|{\Phi}^2_{\varepsilon}(x)-\Phi^2(x)|\leq C\varepsilon;\\
|{f}_{\varepsilon}(x,y_1,y_2,z)-f(x,y_1,y_2,z)|
\leq C\varepsilon;\\
|{g}_{\varepsilon}(x,y_1,y_2,z)-g(x,y_1,y_2,z)|
\leq C\varepsilon,\\
\end{array}
\end{equation*}
$\forall (x,y_1,y_2,z)\in \mathbb{R}\times\mathbb{R}\times\mathbb{R}\times\mathbb{R}$.
\end{lemma}
\textbf{Proof}: We just prove the first inequality. Similarly, we can obtain the other inequalities. By equation (\ref{unae-1}), for fixed $x$, we have,
\begin{equation}
\begin{array}
[c]{ll}
&|{b}_{\varepsilon}(x)-b(x)|\\
\leq & \int_{\mathbb{R}}\left|b(x^0)-b(x)
\right|\frac{1}{\sqrt{2\pi \varepsilon^2}}e^{-\frac{(x^0-x)^2}{2\varepsilon^2}}dx^0\\
\leq & C \int_{\mathbb{R}}\left|x^0-x
\right|\frac{1}{\sqrt{2\pi \varepsilon^2}}e^{-\frac{(x^0-x)^2}{2\varepsilon^2}}dx^0\\
=&C\varepsilon\int_{\mathbb{R}}{|\tilde{x}^0|}\frac{1}{\sqrt{2\pi}}
e^{-\frac{(\tilde{x}^0)^2}{2}}d\tilde{x}^0\\
\leq & C \varepsilon
\end{array}
\end{equation}
where $C$ will change line by line.

This completes the proof. $\ \ \ \ \ \ \ \ \  \Box$

\begin{remark}
\label{re-1}
Notice that $\psi\in\mathcal{C}^{\infty}$, by the property of convolution, we obtain that ${b}_{\varepsilon},{\sigma}_{\varepsilon},\Phi^1_{\varepsilon},\Phi^1_{\varepsilon}
$ satisfy Assumption \ref{ass1}, ${f}_{\varepsilon},{g}_{\varepsilon}$ satisfy Assumption \ref{ass2}, thus we can use the results of \cite{PP92}.
\end{remark}

Now, we show the uniqueness results for the viscosity solution of equation (\ref{vis-1}).
\begin{theorem}
\label{vis-un} Let $b,\sigma$ satisfy  Assumptions \ref{ass-1}; $f,g$ satisfy  Assumptions \ref{ass-2}, \ref{ass-3} and \ref{ass-4}. Then $(u,v)$ is the unique viscosity
solution of the  semi-linear PDEs (\ref{vis-1}).
\end{theorem}
\textbf{Proof}: By Theorem \ref{viscosity}, we obtain the following step:

\textbf{Step 1}. The solution of BSDEs (\ref{pde-1}) and (\ref{pde-2}) $(u,v)$ is the  viscosity
solution of the semi-linear PDEs (\ref{vis-1}).

\textbf{Step 2}. We now prove that the solution  $(u,v)$ is the maximum viscosity subsolution and minimum viscosity supersolution of PDEs (\ref{vis-1}).

Let us consider the following forward backward systems, for $t\in\lbrack0,T]$ and $x\in\mathbb{R}$,
\begin{equation*}
\begin{array}
[c]{ll}
dX_{\varepsilon,s}^{t,x} ={b}_{\varepsilon}(X_{\varepsilon,s}^{t,x})ds
+{\sigma}_{\varepsilon}(X_{\varepsilon,s}^{t,x})dW_{s};\\
dY_{\varepsilon,s}^{1,t,x}= -{f}_{\varepsilon}(X_{\varepsilon,s}^{t,x},Y_{\varepsilon,s}%
^{t,x},Z_{\varepsilon,s}^{1,t,x})ds+Z_{\varepsilon,s}^{1,t,x}dW_{s};\\
dY_{\varepsilon,s}^{2,t,x}=  -{g}_{\varepsilon}(X_{\varepsilon,s}^{t,x},Y_{\varepsilon,s}%
^{t,x},Z_{\varepsilon,s}^{2,t,x})ds+Z_{\varepsilon,s}^{2,t,x}dW_{s};\\
 X_{\varepsilon,t}^{t,x}=x,\ Y_{\varepsilon,T}^{1,t,x}=\Phi^{1}_{\varepsilon}(X_{\varepsilon,T}^{t,x}),\  Y_{\varepsilon,T}^{2,t,x}=\Phi^{2}_{\varepsilon}(X_{\varepsilon,T}^{t,x}),
\end{array}
\end{equation*}
where $Y_{\varepsilon}^{t,x}=(Y_{\varepsilon}^{1,t,x},Y_{\varepsilon}^{2,t,x})$. By Remark \ref{re-1}, ${b}_{\varepsilon},{\sigma}_{\varepsilon},
\Phi^1_{\varepsilon},\Phi^1_{\varepsilon}$ satisfy Assumption \ref{ass1}, ${f}_{\varepsilon},{g}_{\varepsilon}$ satisfy Assumption \ref{ass2}. For given initial data $(t,x)$, we define
\[
u_{\varepsilon}(t,x):=Y_{{\varepsilon},t}^{1,t,x},\ v_{\varepsilon}(t,x):=
Y_{{\varepsilon},t}^{2,t,x}.
\]
Then, by Theorem \ref{2.3}, we obtain $u_{\varepsilon}(t,x), u_{\varepsilon}(t,x)\in \mathcal{C}^{1,2}([0,T]\times\mathbb{R})$ is the  classical solution of the following PDEs:%
\begin{equation}%
\begin{array}
[c]{l}%
\partial_{t}u_{\varepsilon}(t,x)+\mathcal{L}_{\varepsilon}u_{\varepsilon}(t,x)
+f_{\varepsilon}(x,u_{\varepsilon},v_{\varepsilon},\sigma_{\varepsilon}(x)\partial_{x}u_{\varepsilon}(t,x))=0;\\
u_{\varepsilon}(T,x)=\Phi_{\varepsilon}^{1}(x),\text{ \ \ }x\in\mathbb{R};\text{\ }\\
\partial_{t}v_{\varepsilon}(t,x)+\mathcal{L}_{\varepsilon}v_{\varepsilon}(t,x)
+g_{\varepsilon}(x,u_{\varepsilon},v_{\varepsilon},\sigma_{\varepsilon}(x)\partial_{x}v_{\varepsilon}(t,x))=0;\\
v_{\varepsilon}(T,x)=\Phi_{\varepsilon}^{2}(x),\text{ \ \ }x\in\mathbb{R},
\end{array}
\label{vis-2}%
\end{equation}
where
\[
\mathcal{L}_{\varepsilon}=\frac{1}{2}(\sigma_{\varepsilon})^2\partial_{xx}
+b_{\varepsilon}(x)\partial_{x},
\]
and
\[
u_{\varepsilon}(s,X_{{\varepsilon},s}^{t,x})=Y_{{\varepsilon},s}^{1,t,x},\ v_{\varepsilon}(s,X_{{\varepsilon},s}^{t,x})=Y_{{\varepsilon},s}^{2,t,x}.
\]
By Lemma \ref{visle-1} and Assumptions \ref{ass-1} and \ref{ass-2}, there is a constant $C$ such that
\begin{equation}%
\begin{array}
[c]{l}%
\partial_{t}u_{\varepsilon}(t,x)+\mathcal{L}_{}u_{\varepsilon}(t,x)
+f_{}(x,u_{\varepsilon},v_{\varepsilon},\sigma_{}(x)\partial_{x}u_{\varepsilon}(t,x))\leq C\varepsilon;\\
\partial_{t}v_{\varepsilon}(t,x)+\mathcal{L}_{}v_{\varepsilon}(t,x)
+g_{}(x,u_{\varepsilon},v_{\varepsilon},\sigma_{}(x)\partial_{x}v_{\varepsilon}(t,x))\leq C\varepsilon,\\
\end{array}
\label{vis-3}%
\end{equation}
and
\begin{equation}
\label{nede-0}
|u(t,x)-u_{\varepsilon}(t,x)|+|v(t,x)-v_{\varepsilon}(t,x)|
\leq C\varepsilon.
\end{equation}
Then, by inequality (\ref{nede-0}), we obtain
\begin{equation}%
\begin{array}
[c]{l}%
\partial_{t}{u}_{\varepsilon}(t,x)+\mathcal{L}_{}{u}_{\varepsilon}(t,x)
+f_{}(x,{u}_{},{v}_{},\sigma_{}(x)\partial_{x}{u}_{\varepsilon}(t,x))\leq C\varepsilon;\\
\partial_{t}{v}_{\varepsilon}(t,x)+\mathcal{L}_{}{v}_{\varepsilon}(t,x)
+g_{}(x,{u}_{},{v}_{},\sigma_{}(x)\partial_{x}{v}_{\varepsilon}(t,x))\leq C\varepsilon.\\
\end{array}
\label{vis-4}%
\end{equation}
Now, let us assume that $(u_1,v_1)$ is one of the viscosity subsolution of PDEs (\ref{vis-1}), and recall that $(u,v)$ is the viscosity solution of PDEs (\ref{vis-1}). Since ${u}_{\varepsilon},{u}_{1},{v}_{\varepsilon},{v}_{1}$ are continuous functions in $[0,T]\times \mathbb{R}$, we suppose that ${u}_{\varepsilon}+\alpha q_1(t,x)-{u}_{1}$ and ${v}_{\varepsilon}+\alpha q_2(t,x)-{v}_{1}$  take  minimum values  $K_1^{\alpha}$ and $K_2^{\alpha}$ in   $[0,T]\times \mathbb{R}$ at $(t_1,x_1)$ and $(t_2,x_2)$, where $0<\alpha$ is a constant, and $0\leq q_1,q_2 \in \mathcal{C}^{1,2}([0,T]\times\mathbb{R})$.

If $K_1^{\alpha},K_2^{\alpha}\geq0$, by inequality (\ref{nede-0}), we have for any $(t,x)\in[0,T]\times\mathbb{R}$
\begin{equation}
\begin{array}
[c]{ll}
&{u}_{}(t,x)+\alpha q_1(t,x)-{u}_{1}(t,x)\\
\geq&{u}_{\varepsilon}(t,x)+\alpha q_1(t,x)-{u}_{1}(t,x)+C\varepsilon\\
\geq& {u}_{\varepsilon}(t_1,x_1)+\alpha q_1(t_1,x_1)-{u}_{1}(t_1,x_1)+C\varepsilon\\
=&K_1^{\alpha}+C\varepsilon\\
\geq&C\varepsilon.
\end{array}
\end{equation}
Letting $\varepsilon,\alpha \to 0$, thus, ${u}_{}(t,x)-{u}_{1}(t,x)\geq 0$, similarly we have ${v}_{}(t,x)-{v}_{1}(t,x)\geq 0$.

It is easy to verify that $K_1^{\alpha},K_2^{\alpha}$ are decreasing in $\alpha$. We next show that the limitation of $K_1^{\alpha},K_2^{\alpha}$ exist and larger than $0$,  denote
\begin{equation}
\label{nede-2}
\bar{u}_{\varepsilon}(t,x)={u}_{\varepsilon}(t,x)+\alpha q_1(t,x)-K_1^{\alpha},\
\bar{v}_{\varepsilon}(t,x)={v}_{\varepsilon}(t,x)+\alpha q_2(t,x)-K_2^{\alpha}.
\end{equation}
Therefore, by inequality (\ref{vis-4}) we have the following equality about $\bar{u}_{\varepsilon}(t,x),\bar{v}_{\varepsilon}(t,x)$,
\begin{equation}%
\begin{array}
[c]{ll}%
\partial_{t}\bar{u}_{\varepsilon}(t_1,x_1)+\mathcal{L}_{}\bar{u}_{\varepsilon}(t_1,x_1)
+f_{}(x_1,{u}_{},{v}_{},\sigma_{}(x_1)\partial_{x}\bar{u}_{\varepsilon}(t_1,x_1))
\leq C(\varepsilon+\alpha);\\
\partial_{t}\bar{v}_{\varepsilon}(t_2,x_2)+\mathcal{L}_{}\bar{v}_{\varepsilon}(t_2,x_2)
+g_{}(x_2,{u}_{},{v}_{},\sigma_{}(x_2)\partial_{x}\bar{v}_{\varepsilon}(t_2,x_2)) \leq  C(\varepsilon+\alpha).\\
\end{array}
\label{vis-41}%
\end{equation}
Notice that, $(u_1,v_1)$ is the viscosity subsolution of PDEs (\ref{vis-1}), by  equation (\ref{nede-2}) , we obtain
\begin{equation*}%
\begin{array}
[c]{l}%
\partial_{t}\bar{u}_{\varepsilon}(t_1,x_1)+\mathcal{L}_{}\bar{u}_{\varepsilon}
(t_1,x_1)
+f_{}(x_1,\bar{u}_{\varepsilon},\bar{v}_{\varepsilon},\sigma_{}(x_1)\partial_{x}\bar{u}_{\varepsilon}(t_1,x_1))
\geq 0;\\
\partial_{t}\bar{v}_{\varepsilon}(t_2,x_2)+\mathcal{L}_{}\bar{v}_{\varepsilon}(t_2,x_2)
+g_{}(x_2,\bar{u}_{\varepsilon},\bar{v}_{\varepsilon},\sigma_{}(x_2)\partial_{x}\bar{v}_{\varepsilon}(t_2,x_2))
\geq 0.\\
\end{array}
\end{equation*}
By equation (\ref{nede-2}) and inequality (\ref{nede-0}), there exists  constant $C$ yields
\begin{equation}%
\begin{array}
[c]{ll}%
&\partial_{t}\bar{u}_{\varepsilon}(t_1,x_1)+\mathcal{L}_{}\bar{u}_{\varepsilon}
(t_1,x_1)+f_{}(x_1,{u}_{}-K_1^{\alpha},{v}_{}-K_2^{\alpha},\sigma_{}(x_1)\partial_{x}\bar{u}_{\varepsilon}(t_1,x_1))\geq -C(\varepsilon+\alpha);\\
&\partial_{t}\bar{v}_{\varepsilon}(t_2,x_2)+\mathcal{L}_{}\bar{v}_{\varepsilon}
(t_2,x_2)+g_{}(x_2,{u}_{}-K_1^{\alpha},{v}_{}-K_2^{\alpha},\sigma_{}(x_2)\partial_{x}\bar{v}_{\varepsilon}(t_2,x_2))
\geq -C(\varepsilon+\alpha).\\
\end{array}
\label{vis-5}
\end{equation}
 Then, by Assumption \ref{ass-2}, $f,g$ are Lipschatiz continuous functions. Combining  equations (\ref{vis-41}) and  (\ref{vis-5}), we have
\begin{equation}%
\label{vis-6}
\begin{array}
[c]{l}
f_{}(x_1,{u}_{}-K_1^{\alpha},{v}_{}-K_2^{\alpha},\sigma_{}(x_1)\partial_{x}\bar{u}_{\varepsilon}(t_1,x_1)) \geq f_{}(x_1,{u}_{},{v}_{},\sigma_{}(x_1)\partial_{x}\bar{u}_{\varepsilon}(t_1,x_1)) -2C(\varepsilon+\alpha);\\
g_{}(x_2,{u}_{}-K_1^{\alpha},{v}_{}-K_2^{\alpha},\sigma_{}(x_2)\partial_{x}\bar{u}_{\varepsilon}(t_2,x_2)) \geq g_{}(x_2,{u}_{},{v}_{},\sigma_{}(x_2)\partial_{x}\bar{u}_{\varepsilon}(t_2,x_2)) -2C(\varepsilon+\alpha).
\end{array}
\end{equation}
By Assumptions \ref{ass-2} and \ref{ass-4}, there exists Lipschatiz constants $C_2>C_1>0$ such that
\begin{equation*}%
\begin{array}
[c]{l}
f_{}(x_1,{u}_{}-K_1^{\alpha},{v}_{}-K_2^{\alpha},\sigma_{}(x_1)\partial_{x}\bar{u}_{\varepsilon}(t_1,x_1)) - f_{}(x_1,{u}_{},{v}_{},\sigma_{}(x_1)\partial_{x}\bar{u}_{\varepsilon}(t_1,x_1)) \leq C_2K_1^{\alpha}-C_1K_2^{\alpha};\\
g_{}(x_2,{u}_{}-K_1^{\alpha},{v}_{}-K_2^{\alpha},\sigma_{}(x_2)\partial_{x}\bar{u}_{\varepsilon}(t_2,x_2)) -g_{}(x_2,{u}_{},{v}_{},\sigma_{}(x_2)\partial_{x}\bar{u}_{\varepsilon}(t_2,x_2)) \leq C_2K_2^{\alpha}-C_1K_1^{\alpha}.\end{array}
\end{equation*}
which deduce that
\begin{equation*}%
\begin{array}
[c]{l}
-2C(\varepsilon+\alpha)\leq C_2K_1^{\alpha}-C_1K_2^{\alpha};\\
-2C(\varepsilon+\alpha)\leq C_2K_2^{\alpha}-C_1K_1^{\alpha},\end{array}
\end{equation*}
Thus,
$$
-4C(\varepsilon+\alpha)\leq (C_2-C_1)(K_1^{\alpha}+K_2^{\alpha}).\\
$$
Let $\varepsilon,\alpha \to 0$,  notice that $K_1^{\alpha},K_2^{\alpha}$ are decreasing in $\alpha$, thus, the limitation of $K_1^{\alpha},K_2^{\alpha}$ exist, and we obtain
$$
0\leq K_1+K_2.\\
$$
where $K_1=\lim_{\alpha\to 0}K_1^{\alpha}$, $K_2=\lim_{\alpha\to 0}K_2^{\alpha}$.

Now, let us assume $K_1\cdot K_2\leq0$, if $K_1=K_2=0$, then we obatin the assertion. Thus, we only need consider the case: $K_1\geq 0$, $K_2<0$ or $K_2\geq 0$, $K_1<0$, let $\varepsilon,\alpha$ small enough, by Assumption \ref{ass-4}, we have if $K_1\geq 0$, $K_2<0$
$$
g_{}(x_2,{u}_{}-K_1^{\alpha},{v}_{}-K_2^{\alpha},\sigma_{}(x_2)\partial_{x}\bar{u}_{\varepsilon}(t_2,x_2)) -g_{}(x_2,{u}_{},{v}_{},\sigma_{}(x_2)\partial_{x}\bar{u}_{\varepsilon}(t_2,x_2))<0,
$$
similarly, if $K_2\geq 0$, $K_1<0$, thus
$$
f_{}(x_1,{u}_{}-K_1^{\alpha},{v}_{}-K_2^{\alpha},\sigma_{}(x_1)\partial_{x}\bar{u}_{\varepsilon}(t_1,x_1)) - f_{}(x_1,{u}_{},{v}_{},\sigma_{}(x_1)\partial_{x}\bar{u}_{\varepsilon}(t_1,x_1))<0,
$$
which is contrary to the inequality (\ref{vis-6}). Therefore, we obtain that for any $(t,x)\in [0,T]\times \mathbb{R}$,
$$
{u}_{1}(t,x)\leq{u}_{}(t,x),\ {v}_{1}(t,x)\leq{v}_{}(t,x).
$$

Thus, we prove that $(u,v)$ is the maximum viscosity subsolution of PDEs (\ref{vis-1}). Using the similar method as above, we can prove that $(u,v)$ is the minimum viscosity supersolution of PDEs (\ref{vis-1}).

\textbf{Step 3}. Let us assume that $(u_2,v_2)$  is the viscosity solution of PDEs (\ref{vis-1}). Thus, $(u_2,v_2)$  is the viscosity subsolution  and supersolution of PDEs (\ref{vis-1}), following the results of Step 2, we have
$$
u_2(t,x)\leq u(t,x) \leq u_2(t,x),\ v_2(t,x)\leq v(t,x) \leq v_2(t,x),
$$
thus,  $u_2=u,\ v_2=v$, for any $(t,x)\in [0,T]\times\mathbb{R}$.

This completes the proof. $\ \ \ \ \  \ \Box$

\begin{remark}
Now, let us extend our model as follows:
\begin{align}
dX_{s}^{t,x}  &  =b(s,X_{s}^{t,x})ds+\sigma(s,X_{s}^{t,x})dW_{s}%
,\label{rpde-1}\\
X_{t}^{t,x}  &  =x,\nonumber
\end{align}
and%
\begin{equation}%
\begin{array}
[c]{rl}%
dY_{s}^{1,t,x}= &-f_1(s,X_{s}^{t,x},Y_{s}%
^{t,x},Z_{s}^{1,t,x})ds+Z_{s}^{1,t,x}dW_{s};\\
\cdots & \cdots\\
dY_{s}^{i,t,x}=& -f_i(s,X_{s}^{t,x},Y_{s}%
^{t,x},Z_{s}^{i,t,x})ds+Z_{s}^{i,t,x}dW_{s};\\
\cdots & \cdots\\
dY_{s}^{n,t,x}=& -f_n(s,X_{s}^{t,x},Y_{s}%
^{t,x},Z_{s}^{n,t,x})ds+Z_{s}^{n,t,x}dW_{s};\\
Y_{T}^{i,t,x}=&\Phi^{i}(X_{T}^{t,x}),\ i=1,2,\cdots,n.
\end{array}
\label{rpde-2}%
\end{equation}
where $Y^{t,x}=(Y^{1,t,x},\cdots,Y^{i,t,x},\cdots,Y^{n,t,x})$ and
\[%
\begin{array}
[c]{l}%
b, \sigma:[0,T]\times\mathbb{R}^n\rightarrow\mathbb{R}^n;\\
f_i:[0,T]\times\mathbb{R}^n\times\mathbb{R}^n\times\mathbb{R}\rightarrow\mathbb{R}^{};\\
\Phi^{i}:\mathbb{R}^n\rightarrow\mathbb{R},
\end{array}
\]
define that
\[
u_i(t,x)=Y_t^{i,t,x},\ i=1,2,\cdots,n.
\]
Thus, the related semi-linear PDEs are
\begin{equation}%
\begin{array}
[c]{l}%
\partial_{t}u_i(t,x)+\mathcal{L}u_i(t,x)+f_i(x,u,\sigma^T(x)\partial_{x}u_i(t,x))=0;\\
u_i(T,x)=\Phi^{i}(x),\text{ \ \ }x\in\mathbb{R},\ i=1,2,\cdots,n,\text{\ }\\
\end{array}
\label{rvis-1}%
\end{equation}
where $u=(u_1,u_2,\cdots,u_n)$ and
\[
\mathcal{L}=\frac{1}{2}\sigma\sigma^{T}(t,x)\partial_{xx}+b(t,x)\partial_{x}.
\]
Notice that, if we suppose that $b,\sigma,f_i,\Phi^i$ satisfy the same conditions as in this study, we can prove that $u$ is the unique viscosity solution of PDEs (\ref{rvis-1}).
\end{remark}

\appendix
\section{The proof of Theorem \ref{theorem-cp}}
\textbf{Proof}: By the basic theory of BSDEs with coefficients are Lipschtiz continuous and linear growth in $(y,z)$, equations (\ref{cp-11}) and (\ref{cp-12}) have unique solution, we refer \cite{EPQ97} and \cite{PP90} for basic theory of BSDEs.

Therefore, we denote the solutions of equations (\ref{cp-11}%
) and (\ref{cp-12}) as $(Y_{i},Z_{i}).$ Setting $\hat{Y}_i={Y}_{i,2}-{Y}_{i,1},$
$\hat{Z}_i=Z_{i,2}-Z_{i,1},$ $\hat{f}(\cdot)=f_{2}-f_{1},$
$\hat{\xi}_{i}=\xi_{i,2}-\xi_{i,1},i=1,2;\ \hat{Z}=(\hat{Z}_1,\hat{Z}_2).$

Firstly, we consider equation (\ref{cp-11}), $(\hat{Y}_1,\hat{Z}_{1})$ satisfies the
following equation%
\[%
\begin{array}
[c]{rl}%
\hat{Y}_{1}(t)= & \hat{\xi}_1+\int_{t}^{T}[\hat{f}(Y_{1}(s),Z_{1,1}(s))+f_2(Y_{1,1}(s),Y_{2,2}(s)
,Z_{1,1}(s)) -f_2(Y_{1,1}(s),Y_{2,1}(s),Z_{1,1}(s))\\
&+a(s)\hat{Y}_{1}(s)+b(s)\hat{Z}_{1}(s)]ds-\int%
_{t}^{T}\hat{Z}_{1}(s)d{W}(s),
\end{array}
\]
where
\begin{equation*}%
\begin{array}
[c]{rl}%
a(s)=
\genfrac{\{}{.}{0pt}{}{\frac{f_{2}({Y}_{1,2}(s))-f_{2}%
({Y}_{1,1}(s))}{{Y}_{1,2}(s)-{Y}%
_{1,1}(s)},\text{ \ }Y_{1,2}(s)-Y_{1,1}(s)\neq
0}{\text{\ \ \ \ \ \ \ }0,\text{  \ \ \ \ \  \ \ \ \ \  \ \  \ \ \ \ \ \ \ \ \ \ }%
{Y}_{1,2}(s)-{Y}_{1,1}(s)=0};\\
b(s)=
\genfrac{\{}{.}{0pt}{}{\frac{f_{2}({Z}_{1,2}(s))-f_{2}%
({Z}_{1,1}(s))}{{Z}_{1,2}(s)-{Z}%
_{1,1}(s)},\text{ \ }Z_{1,2}(s)-Z_{1,1}(s)\neq
0}{\text{\ \ \ \ \ \ \ }0,\text{  \ \ \ \ \  \ \ \ \ \ \ \ \ \ \ \ \ \ \ \ \ }%
{Z}_{1,2}(s)-{Z}_{1,1}(s)=0}.
\end{array}
\end{equation*}
Notice that $f_{2}$ satisfies Lipschitz condition and $f_2\geq f_1$,
thus $\left\vert a(s)\right\vert + \left\vert b(s)\right\vert\leq 2C$ and $\hat{f}(Y_{1}(s),Z_{1,1}(s))=f_2(Y_{1,1}(s),Y_{2,1}(s),Z_{1,1}(s))-f_1(Y_{1,1}(s),
Y_{2,1}(s),Z_{1,1}(s))\geq0.$

Consider the following SDE:
\[%
\begin{array}
[c]{ll}%
dX(s)=-X(s)a(s)ds+X(s)b(s)d{W}(s);\\
X(0)=  1.
\end{array}
\]
Applying It\^{o} formula to $\hat{Y}_1(t)X(t)$, we have
\begin{equation}
\begin{array}
[c]{rl}%
\hat{Y}_{1}(t)X(t)-\hat{\xi}_{1}{X(T)}{}
=&\int_{t}^{T}[\hat{f}(Y_{1}(s),Z_{1,1}(s))
+f_2(Y_{1,1}(s),Y_{2,2}(s),Z_{1,1}(s))\\
&-f_2(Y_{1,1}(s),Y_{2,1}(s),Z_{1,1}(s))]{X(s)}ds\\
&-\int_t^TZ_{1,1}(s)X(s)dW(s)
+\int_t^T\hat{Y}_{1}(s)X(s)b(s)dW(s).\\
\label{cp-50}
\end{array}
\end{equation}
Taking conditional expectation on both side of equation (\ref{cp-50}), we have
\begin{equation*}
\begin{array}
[c]{rl}%
&({Y}_{1,2}(t)-Y_{1,1}(t))X(t)-E[\int_{t}^{T}[f_2(Y_{1,1}(s),Y_{2,2}(s),Z_{1,1}(s))
-f_2(Y_{1,1}(s),Y_{2,1}(s),Z_{1,1}(s))]
{X(s)}{}ds\mid \mathcal{F}_t]\\
=&E[\hat{\xi}_1X(T)+\int_t^T\hat{f}(Y_{1}(s),Z_{1,1})(s)X(s)ds\mid \mathcal{F}_t],
\end{array}
\end{equation*}
which deduce that
\begin{equation*}
\begin{array}
[c]{rl}%
&({Y}_{1,2}(t)-Y_{1,1}(t))X(t)-E[\int_{t}^{T}[f_2(Y_{1,1}(s),Y_{2,2}(s),Z_{1,1}(s))
-f_2(Y_{1,1}(s),Y_{2,1}(s),Z_{1,1}(s))]
{X(s)}{}ds\mid \mathcal{F}_t]\geq 0.
\end{array}
\end{equation*}

 For convenience, we denote
 $$
 f_2(Y_{2,2}(s))=f_2(Y_{1,1}(s),Y_{2,2}(s),Z_{1,1}(s));
 \quad f_2(Y_{2,1}(s))=f_2(Y_{1,1}(s),Y_{2,1}(s),Z_{1,1}(s)),
 $$
 thus
\begin{equation}
{Y}_{1,2}(t)-Y_{1,1}(t)\geq E[\int_{t}^{T}[f_2(Y_{2,2}(s))-f_2(Y_{2,1}(s))]
\frac{X(s)}{X(t)}ds\mid \mathcal{F}_t].
\label{cpp-1}
\end{equation}
Using the same method as above, we have
\begin{equation}
{Y}_{2,2}(t)-Y_{2,1}(t)\geq E[\int_{t}^{T}[g_2(Y_{1,2}(s))-g_2(Y_{1,1}(s))]
\frac{\bar{X}(s)}{\bar{X}(t)}ds\mid \mathcal{F}_t],
\label{cpp-2}
\end{equation}
where
\begin{equation*}
\begin{array}
[c]{ll}
g_2(Y_{1,2}(s))=g_2(Y_{1,2}(s),Y_{2,2}(s),Z_{2,1}(s));\\
g_2(Y_{1,1}(s))=g_2(Y_{1,1}(s),Y_{2,2}(s),Z_{2,1}(s)),
\end{array}
\end{equation*}
the $\bar{X}(s)$ is the exponential martingale similar with $X(s)$, $0\leq s\leq T$.

From equations (\ref{cpp-1}) and (\ref{cpp-2}), we have $Y_{1,2}(T)\geq Y_{1,1}(T)$ and $Y_{2,2}(T)\geq Y_{2,1}(T)$. By Assumption \ref{ass-3}, $f_2$ is non-decreasing in the second dimension of argument, $g_2$ is non-decreasing in the first dimension of argument, by Gronwall inequality, we get
$$
{Y}_{1,2}(t)\geq Y_{1,1}(t),\quad {Y}_{2,2}(t)\geq Y_{2,1}(t).
$$

This completes the proof. $\ \ \ \ \  \ \Box$

\end{document}